\documentclass[12pt,oneside]{amsart}
\usepackage{geometry}                
\usepackage{graphicx}
\usepackage{amssymb}
\usepackage{epstopdf}
\usepackage{amsmath,amsthm,amscd,amssymb}
\usepackage{latexsym}
\usepackage[colorlinks,citecolor=red,pagebackref,hypertexnames=false]{hyperref}

\numberwithin{equation}{section}

\theoremstyle{plain}
\newtheorem{theorem}{Theorem}[section]
\newtheorem{lemma}[theorem]{Lemma}

\theoremstyle{definition}
\newtheorem{definition}[theorem]{Definition}

\theoremstyle{remark}
\newtheorem{remark}[theorem]{Remark}

\newtheorem{case[theorem]}{Case}

\title{On directions determined by subsets of vector spaces over finite fields}
\date{October 4, 2010}      
\author{Alex Iosevich, Hannah Morgan and Jonathan Pakianathan}
\address{Department of Mathematics, University of Rochester, Rochester, NY}
\email{iosevich@math.rochester.edu}
\address{Department of Mathematocs, University of Missouri, Columbia, MO}
\email{hmm5x7@mizzou.edu}
\address{Department of Mathematics, University of Rochester, Rochester, NY} 
\email{jonpak@math.rochester.edu} 

\thanks{This work was partially supported by the NSF Grant DMS10-45404. The article is based in part on Hannah Morgan's senior thesis at the University of Missouri-Columbia}

\thanks{The authors wish to thank Doowon Koh, Misha Rudnev, Steven Senger, Chun-Yen Shen and the anonymous referee for some very helpful remarks and suggestions} 

\dedicatory{This paper is dedicated to the memory of Nigel Kalton} 

\begin{document}
\maketitle

\begin{abstract} We prove that if a subset of a $d$-dimensional vector space over a finite field with $q$ elements has more than $q^{d-1}$ elements, then it determines all the possible directions. We obtain a complete characterization if the size of the set is $\ge q^{d-1}$. If a set has more than $q^k$ elements, its direction set projects onto a $k$-dimensional set. We prove stronger results for sets that are sufficiently random. This result is best possible as the example of a $k$-dimensional hyperplane shows. We can view this question as an Erd\H os type problem where a sufficiently large subset of a vector space determines a large number of configurations of a given type. See, for example, (\cite{PS04}), \cite{BKT04}, \cite{HIKR07}, \cite{TV06}, \cite{V08} and the references contained therein. For discrete subsets of ${\Bbb R}^d$, this question has been previously studied by Pach, Pinchasi and Sharir. See (\cite{PPS07}). \end{abstract}                

\section{Introduction} 

The celebrated Kakeya conjecture, proved in the finite field context by Dvir (\cite{D09}), says that if $E \subset {\Bbb F}_q^d$, $d \ge 2$, contains a line (or a fixed positive proportion thereof) in every possible direction, then $|E| \ge cq^d$. Here, and throughout, $|E|$ denotes the number of elements of $E$ and ${\Bbb F}_q^d$ denotes the $d$-dimensional vector space over the finite field with $q$ elements. 

While Dvir's theorem shows that a set containing a line in every directions is large, in this paper we see to determine how large a set needs to be to determine every possible direction, or a positive proportion thereof. In the discrete setting the problem of directions was studied in recent years by Pach, Pinchasi, and Sharir. See \cite{PS04} and \cite{PPS07}. In the latter paper they prove that if $P$ is a set of $n$ points in ${\Bbb R}^3$, not all in a common plane, then the pairs of points of $P$ determine at least $2n-5$ distinct directions if $n$ is odd and at least $2n-7$ distinct directions if $n$ is even. 

In order to state our main result, we need to make precise the notion of directions in subsets of ${\Bbb F}_q^d$. 

\begin{definition} We say that two vectors $x$ and $x'$ in ${\Bbb F}_q^d$ point in the same direction if there exists $t \in {\Bbb F}_q^{*}$ such that $x'=tx$. Here ${\Bbb F}_q^{*}$ denotes the multiplicative group of ${\Bbb F}_q$. Writing this equivalence as $x \sim x'$, we define the set of directions as the quotient 
\begin{equation} \label{directionequivalence} {\mathcal D}({\Bbb F}_q^d)={\Bbb F}_q^d / \sim. \end{equation}

Similarly, we can define the set of directions determined by $E \subset {\Bbb F}_q^d$ by 

\begin{equation} \label{equivalence} {\mathcal D}(E)=E-E / \sim, \end{equation} where 
$$ E-E=\{x-y: x,y \in E\},$$ with the same equivalence relation $\sim$ as in (\ref{directionequivalence}) above. 
\end{definition} 

It is not difficult to see that $|{\mathcal D}({\Bbb F}_q^d)|=q^{d-1}(1+o(1))$. Thus the question above is rephrased in the following form. How large does $E \subset {\Bbb F}_q^d$ need to be to ensure that ${\mathcal D}(E)={\mathcal D}({\Bbb F}_q^d)$, or, more modestly, that ${\mathcal D}(H_n) \subset {\mathcal D}(E)$, where $H_n$ is a $n$-dimensional plane. 

Since $E$ may be a $k$-dimensional plane, a {\it necessary} condition for ${\mathcal D}(H_{k+1}) \subset {\mathcal D}(E)$ is $|E|>q^{k}$. We shall see that this simple necessary condition is in fact sufficient. Our main result is the following. 
\begin{theorem} \label{main} Let $E \subset {\Bbb F}_q^d$.  Suppose that $|E|>q^k$, $1 \leq k \leq d-1$. Then there exists a $(k+1)$-dimensional subspace, denoted by $H_{k+1}$ such that ${\mathcal D}(E)$ projects onto $H_{k+1}$. In particular, if $|E|>q^{d-1}$, every possible direction is determined. \end{theorem}

It is reasonable to conjecture that if $|E|=q^k$, then $|{\mathcal D}(E)| \gtrsim q^k$ \footnote{Here and throughout, $X \lesssim Y$ means that there exists $C>0$, independent of $q$, such that $X \leq CY$.} unless $E$ is of a certain characterized form. In the case $k=d-1$, we have the following characterization. 

\begin{theorem} \label{graph} Suppose that $E \subset {\Bbb F}_q^d$ has ${\mathcal D}(E) \not={\mathcal D}({\Bbb F}_q^d)$. Then after a linear change of coordinates, we have 
\begin{equation} \label{characterization} E=\{(x_1, \dots, x_{d-1}, f(x_1, \dots, x_{d-1})): (x_1, \dots, x_{d-1}) \subset {\mathcal U} \} \end{equation} for some $f: {\mathcal U} \to {\Bbb F}_q$, where ${\mathcal U} \subset {\Bbb F}_q^{d-1}$. 

Conversely, any $E$ of this form has the property that ${\mathcal D}(E) \not={\mathcal D}({\Bbb F}_q^d)$. 
\end{theorem} 

Observe that the size of $E$ in (\ref{characterization}) is at most $q^{d-1}$. Therefore, we have a complete characterization of the situation when $|E| \ge q^{d-1}$. 

It would be very interesting to prove a version of this result where instead of considering directions determined by a single set $E$, we consider directions determined by pairs of points where one lies in $E \subset {\Bbb F}_q^d$ and the other in $F \subset {\Bbb F}_q^d$. We are able to do this in the case when $k=d-1$. 

\begin{theorem} \label{differentsets} Let $E, F \subset {\Bbb F}_q^d$, not necessarily disjoint. Define ${\mathcal D}(E,F)=E-F / \sim$, where $\sim$ is as in (\ref{equivalence}) above. Suppose that $|E|+|F|>q^d$ {\bf or} $|E \cap F| > q^{d-1}$, then ${\mathcal D}(E,F)={\mathcal D}({\Bbb F}_q^d)$. Moreover, the result is sharp in the sense that there exist $E,F \subset {\Bbb F}_q^d$ with $|E|+|F|=q^d$, such that ${\mathcal D}(E,F) \not={\mathcal D}({\Bbb F}_q^d)$. 
\end{theorem} 

\begin{remark} The proof of Theorem \ref{differentsets} below actually characterizes the sets in question, but the characterization is not particularly pretty. \end{remark} 

Theorem \ref{main} is in general best possible as we note above. However, if the set is sufficiently "random", we can obtain stronger conclusions. One reasonable measure of randomness of a set is via the size of its Fourier coefficients. Let $\chi$ denote a non-trivial principal character on ${\Bbb F}_q$. See \cite{LN97} for a thorough description of this topic. Note that if $q$ is prime, we can take 
$\chi(t)=e^{\frac{2 \pi i t}{q}}.$ The basic properties of characters are the facts that $\chi(0)=1$, $||\chi(t)||=1$, where $||\cdot||$ denotes complex modulus and 
\begin{equation} \label{orthogonality} q^{-1} \sum_{a \in {\Bbb F}_q} \chi(at)=\delta_0(t),\end{equation} where $\delta_0(t)=1$ if $t=0$ and $0$ otherwise. Given $f: {\Bbb F}_q^d \to {\Bbb C}$, define the Fourier transform of $f$, 
$$ \widehat{f}(m)=q^{-d} \sum_{x \in {\Bbb F}_q^d} \chi(-x \cdot m) f(x).$$ 

We shall also make use of the Plancherel formula 
\begin{equation} \label{plancherel} \sum_m {|\widehat{f}(m)|}^2=q^{-d} \sum_x {|f(x)|}^2. \end{equation}

We are now ready to define the notion of randomness we are going to use. 
\begin{definition} We say that $E \subset {\Bbb F}_q^d$, $d \ge 2$, is a {\it Salem} set if 
\begin{equation} \label{salem} |\widehat{E}(m)| \lesssim q^{-d} \sqrt{|E|} \ \text{for} \ m \not=(0, \dots, 0). \end{equation} 
\end{definition} 

See \cite{IR07} for this definition and examples of Salem sets in ${\Bbb F}_q^d$. See also \cite{W03} where the original version of the concept, in the context of measures in Euclidean space is described. 

Our second result illustrates that under our randomness assumption, we can obtain better exponents. While this is certainly less compelling than the main result above, it does underline the fact that the main obstruction to obtaining many distinct directions, the hyper-plane, has a rather special structure and is far from random. 
\begin{theorem} \label{mainsalem} Suppose that $E \subset {\Bbb F}_q^d$, $d \ge 2$, is a Salem set. 

i) If $|E|>q^{d-1}$, then ${\mathcal D}(E)={\mathcal D}({\Bbb F}_q^d).$ 

ii) If $|E| \leq q^{d-1}$, 
$$ |{\mathcal D}(E)| \gtrsim \min \left\{ \frac{{|E|}^2}{q}, q^{d-1} \right\}.$$ 

iii) If $|E| \leq q^{d-1}$, then 
$$ |{\mathcal D}(E)| \gtrsim |E|.$$ 
\end{theorem} 

In particular, if $|E| \gtrsim q^{\frac{d}{2}}$, $|{\mathcal D}(E)| \gtrsim q^{d-1}$. The lower bound in Part ii) is better than the lower bound in Part iii) if $|E| \gtrsim q$. 

\begin{remark} For a related result, see \cite{HLS10}, Theorem 2.2. See that paper and also \cite{V08} and \cite{TV06} for a connection between the problem under consideration here and the expansion phenomenon in graphs. 
\end{remark} 

\begin{remark} Part ii) holds without the assumption that $|E| \leq q^{d-1}$. We simply wish to emphasize the fact that if $|E|>q^{d-1}$, a much stronger conclusion is already available from Theorem \ref{main}. \end{remark} 

\begin{remark} The proof of Part ii) below does not use the full strength of the Salem assumption (\ref{salem}). What is required is a weaker property $|E-E| \gtrsim \min \{{|E|}^2, q^d \}$, which follows from (\ref{salem}) as Lemma \ref{bigdifferencelemma} shows. To construct a set satisfying this weaker property which is not Salem, just construct a Salem set (see e.g. \cite{IR07}) on ${\Bbb F}_q^k$, $k<d$, and embedd this ${\Bbb F}_q^k$ as a sub-space of ${\Bbb F}_q^d$. 
\end{remark} 

\begin{remark} We note that the conclusion of Theorem \ref{mainsalem} does not in general hold if $E$ is not a Salem set. To see this, take $E \subset H_{k+1}$, $1 \leq k \leq \frac{d}{2}$, a $(k+1)$-dimensional sub-space of ${\Bbb F}_q^d$. Further suppose that $|E| \approx q^{k+\alpha}$ for some $\alpha>0$. Since $E \subset H_{k+1}$, $|{\mathcal D}(E)| \leq q^k$. Therefore, it is not true that $|{\mathcal D}(E)| \gtrsim \frac{{|E|}^2}{q}$ since $q^{2k+2\alpha-1}$ is much greater than $q^k$ when $q$ is large if $k \ge 1$. It is also not true in this case that $|{\mathcal D}(E)| \gtrsim |E|$ since $q^{k+\alpha}$ is much greater than $q^k$ when $q$ is large. \end{remark} 

\begin{remark} What we do not know is to what extent Theorem \ref{mainsalem} can be improved. For example, we do not know of a single Salem subset $E$ of ${\Bbb F}_q^d$ of size $>Cq^{\frac{d-1}{2}}$ for which 
$$|{\mathcal D}(E)|=o(q^{d-1}) \footnote{Recall that $X=o(Y)$ for quantities $X$ and $Y$ depending on the parameter $q$, if $\frac{X}{Y} \to 0$ as $q \to \infty$}.$$

It is reasonable to conjecture that if $E$ is Salem and $|E| \ge Cq^{\frac{d-1}{2}}$, with a sufficiently large constant $C>0$, then ${\mathcal D}(E)={\mathcal D}({\Bbb F}_q^d)$. We do not currently know how to approach this question. \end{remark} 

\vskip.25in 
  
\section{Proof of Theorem \ref{main}}
 
By rotating the coordinates, if necessary, we may define $\nu_E(t_1, ..., t_k)$ by the expression
$$ \left|\left\{(x,y) \in E \times E: x_2-y_2=t_1(x_1-y_1), ..., x_{k+1}-y_{k+1}=t_k(x_1-y_1); 
x \not=y \right\}\right|.$$ 
 
Let $\chi$ denote a non-trivial principal additive character on ${\Bbb F}_q$. It follows from (\ref{orthogonality}) that

$$ \nu_E(t_1,..., t_k)=\sum_{\{(x,y): (x_2-y_2)=t_1(x_1-y_1), ... , (x_{k+1}-y_{k+1})=t_k(x_1-y_1); x \not=y \}} 
E(x)E(y)$$
$$=q^{-k} \sum_{s_1, ..., s_k \in {\Bbb F}_q} \sum_{x \not=y \in {\Bbb F}_q^d}  E(x)E(y) \chi(s_1((x_2-y_2)-t_1(x_1-y_1))) ...\chi(s_k((x_{k+1}-y_{k+1})-t_k(x_1-y_1)))$$
$$ =\frac{|E|(|E|-1)}{q^k}$$
$$-q^{-k} 
\sum_{(s_1, ..., s_k) \not=(0, \dots, 0)} \sum_{x=y \in {\Bbb F}_q^d}  E(x)E(y)\chi(s_1((x_2-y_2)-t_1(x_1-y_1))) ...
\chi(s_k((x_{k+1}-y_{k+1})-t_k(x_1-y_1)))$$
$$+q^{-k} 
\sum_{(s_1, ..., s_k) \not=(0, \dots, 0)} \sum_{x,y \in {\Bbb F}_q^d}  E(x)E(y)\chi(s_1((x_2-y_2)-t_1(x_1-y_1))) ...
\chi(s_k((x_{k+1}-y_{k+1})-t_k(x_1-y_1)))$$
\begin{equation} \label{incidence}=\frac{|E|(|E|-1)}{q^k}-|E| \left(\frac{q^k-1}{q^k}\right)+
R(t_1, ... , t_k). \end{equation}

\begin{lemma} \label{positive} With the notation above, $R(t_1. \dots, t_k) \ge 0$. 
\end{lemma} 

To prove the lemma, we see that by the definition of the Fourier transform, we see that $R(t_1, \dots, t_k)$ equals 

$$ q^{2d-k} \sum_{(s_1, \dots, s_k) \not=(0, \dots, 0)} \widehat{E}(s_1t_1+ ... +s_kt_k,-s_1, ... , -s_k, {\vec 0}) 
\widehat{E}(-s_1t_1- ... -s_kt_k, s_1, ... , s_k, {\vec 0})$$
 $$=q^{2d-k} \sum_{(s_1, \dots, s_k) \not=(0, \dots, 0)} \widehat{E}(s_1t_1+ ... +s_kt_k,-s_1,... , -s_k,{\vec 0}) 
 \overline{\widehat{E}(s_1t_1+ ... +s_kt_k,-s_1,... , -s_k,{\vec 0})}$$
$$=q^{2d-k} \sum_{(s_1, \dots, s_k) \not=(0, \dots, 0)} 
{|\widehat{E}(s_1t_1+ ... +s_kt_k,-s_1,... , -s_k, {\vec 0})|}^2 \ge 0.$$ 

This completes the proof of Lemma \ref{positive}. It follows that 
$$  \nu_E(t_1,..., t_k) \ge \frac{|E|(|E|-1)}{q^k}-|E| \left(\frac{q^k-1}{q^k}\right).$$ 

The right hand side is positive as long as $|E|>q^k$ and this completes the proof. 

\vskip.25in 

\section{Proof of Theorem \ref{mainsalem}} 

\vskip.125in 

Part i) follows instantly from Theorem \ref{main}. To prove Part iii), we observe that by the estimate (\ref{incidence}) above, we have 
$$ \nu_E(t_1, \dots, t_{d-1})=\frac{|E|(|E|-1)}{q^{d-1}}-|E| \left(\frac{q^{d-1}-1}{q^{d-1}}\right)+
R(t_1, ... , t_{d-1}).$$

If $|E| \leq q^{d-1}$, 
$$ \nu_E(t_1, \dots, t_{d-1}) \leq 2|E|+R(t_1, \dots, t_{d-1})$$ 
$$ \lesssim 2|E|+q^{d+1} \cdot q^{d-1} \cdot q^{-2d} |E|=|E|,$$ where the second inequality holds by the Salem property (\ref{salem}). It follows that 
$$ {|E|}^2-|E| \leq \sum_{t_1, \dots, t_{d-1}} \nu_E(t_1, \dots, t_{d-1}) \lesssim |{\mathcal D}(E)| \cdot |E|.$$ 

We conclude that 
$$ |{\mathcal D}(E)| \gtrsim |E|$$ and Part iii) is proved. 

To prove Part ii), we need the following observation. 
\begin{lemma} \label{bigdifferencelemma} \label{difference} Suppose that $E \subset {\Bbb F}_q^d$, $d \ge 2$, is a Salem set. Then 
\begin{equation} \label{bigdifference} |E-E| \gtrsim \min \{{|E|}^2, q^d \}. \end{equation} 
\end{lemma} 

To prove the lemma, define the function $\mu(z)$ by the relation 
\begin{equation} \label{relation} \sum_{z \in {\Bbb F}_q^d} f(z) \mu(z)=\sum_{x,y} f(x-y)E(x)E(y).
\end{equation}

Equivalently, one can set 
$$ \mu(z)=\sum_{x-y=z} E(x)E(y)$$ and check that 
$$ \sum_z f(z) \mu(z)=\sum_z f(z) \sum_{x-y=z} E(x)E(y)=\sum_{x,y} f(x-y)E(x)E(y).$$

Observe that $\mu(z) \not=0$ precisely when $z \in E-E$. Taking $f(z)=\chi(-z \cdot m) q^{-d}$ in (\ref{relation}), we see that 
$$ \widehat{\mu}(m)=q^{-d} \sum_{x,y} \chi((x-y) \cdot m) E(x)E(y)=q^d {|\widehat{E}(m)|}^2.$$ 

Applying (\ref{relation}) once again with $f(z)=1$, we get
$$ \sum \mu(z)={|E|}^2.$$ 

It follows that 
$$ {|E|}^4= {\left( \sum_z \mu(z) \right)}^2 \leq |E-E| \cdot \sum_z \mu^2(z)$$
$$=|E-E| \cdot q^d \cdot \sum_m {|\widehat{\mu}(m)|}^2$$
$$=|E-E| \cdot q^{3d} \cdot \sum_m {|\widehat{E}(m)|}^4$$ 
$$=|E-E| \cdot q^{3d} \cdot q^{-4d} \cdot {|E|}^4$$
$$+|E-E| \cdot q^{3d} \cdot \sum_{m \not=(0, \dots, 0)} {|\widehat{E}(m)|}^4$$
$$ \lesssim |E-E| \cdot q^{-d} \cdot {|E|}^4+
|E-E| \cdot q^{3d} \cdot q^{-2d} |E| \sum_m {|\widehat{E}(m)|}^2$$ 
$$=|E-E| \cdot q^{-d} \cdot {|E|}^4+|E-E| \cdot q^{3d} \cdot q^{-2d} \cdot q^{-d} {|E|}^2$$
$$=|E-E| \cdot q^{-d} \cdot {|E|}^4+|E-E| \cdot {|E|}^2.$$ 

In the fourth line above we used the Salem property (\ref{salem}). In the fifth line, we used the Plancherel formula (\ref{plancherel}). We conclude that 
$$ |E-E| \gtrsim \min \{{|E|}^2, q^d \},$$ as claimed.

We are now ready to complete the proof of Part ii). By definition of ${\mathcal D}(E)$ (\ref{directionequivalence}), at most $q$ points of $E-E$ account for a given element of ${\mathcal D}(E)$. It follows that 
$$ |{\mathcal D}(E)| \gtrsim \frac{{|E|}^2}{q}$$ and the proof of Part ii) is complete. 

\vskip.125in 

\section{Proof of Theorem \ref{graph}} 

First a remark to note that this theorem and its proof hold over any field $\mathbf{k}$ 
but we will assume $\mathbf{k}=\mathbb{F}_q$ throughout as that is the primary field of 
concern for this paper.

We will call a subset $E \subset {\Bbb F}_q^d$ "direction deficient" if 
${\mathcal D}(E) \not={\mathcal D}({\Bbb F}_q^d)$.

If $E$ is direction deficient then ${\mathcal D}(E)$ misses a line $L=\hat{v} /  \sim$ 
where $\hat{v}/ \sim$ denotes the unique line through the origin and nonzero vector $\hat{v}$.

Notice that a linear change of variables of the vector space ${\Bbb F}_q^d$ is given 
by a matrix $\mathbb{A} \in GL_d(\Bbb F_q)$ and that such a change of variables takes 
lines to lines. Futhermore if $\mathbb{A}(E)$ denotes the image of $E$ under this change of variables, 
we have $\mathbb{A}(E)-\mathbb{A}(E)=\mathbb{A}(E-E)$ and so $\mathbb{A}(E)$ is direction 
deficient missing a line $\mathbb{A}(L)$ if and only if $E$ is direction deficient missing a line $L$.

Thus without loss of generality, we can assume our direction deficient set $E$ has 
${\mathcal D}(E)$ not contain the line $L$ which is the $x_d$-axis i.e., the line through the vector 
$(0,0,\dots,0,1)$. Thus $(E-E) \cap L = \{ \hat{0} \}$.

Consider the projection $\pi: \Bbb{F}_q^d \to \Bbb{F}_q^{d-1}$ to the first $(d-1)$-coordinates 
with kernel the line $L$. Let $\pi(E)=U \subseteq \Bbb{F}_q^{d-1}$.

$e_1,e_2 \in E$ have $\pi(e_1)=\pi(e_2)$ if and only if $e_1-e_2 \in ker(\pi)=L$. 
As $(E-E) \cap L = \{ \hat{0} \}$ this happens if and only if $e_1 = e_2$.
Thus $\pi$ restricts to an injective map on $E$ and hence to a bijection 
$\pi |_E: E \to \pi(E)=U \subseteq \Bbb{F}_q^{d-1}$.  Now since $\pi$ was a projection map to the 
first $(d-1)$-coordinates, the inverse map $(\pi |_E)^{-1}: U \to E$ has the form 
$(\pi |_E)^{-1}(x_1,\dots,x_{d-1})=(x_1,\dots,x_{d-1}, f(x_1,\dots,x_{d-1}))$ for a function 
$f: U \to \Bbb{F}_q$ which is the $d$th coordinate function of this inverse.

Since $E=(\pi |_E)^{-1}(U)$ we then find that 
$$
E = \{ (x_1,\dots,x_{d-1},f(x_1,\dots,x_{d-1})) | (x_1,\dots,x_{d-1}) \in U \subseteq \Bbb{F}_q^{d-1} \}
$$
as desired. Thus $E$ (up to linear change of variables) 
is the graph set of a function $f: U \subseteq \Bbb{F}_q^{d-1} \to \Bbb{F}_q$.

Conversely let $E$ be any graph set as above for a function $f: U \subseteq \Bbb{F}_q^{d-1} \to \Bbb{F}_q$. Then we claim ${\mathcal D}(E)$ does not contain the line $L$ through $(0,\dots,0,1)$.

Notice if $e_1=(\hat{u}_1,f(\hat{u}_1))$ and $e_2=(\hat{u}_2,f(\hat{u}_2))$ then 
$e_1 - e_2 = (\hat{u}_1-\hat{u}_2, f(\hat{u}_1)-f(\hat{u}_2))$ never is of the form 
$(0,\dots,0,\text{ nonzero })$ as $\hat{u}_1=\hat{u}_2 \to f(\hat{u}_1)=f(\hat{u}_2)$. 
Thus ${\mathcal D}(E)$ does not contain the line $L$ through $(0,\dots,0,1)$ 
and so we see any graph set is direction deficient and the proof of the theorem is complete.

\vskip.125in 

\section{Proof of Theorem \ref{differentsets}} 

Let $E,F \subseteq \Bbb{F}_q^d$ not necessarily disjoint. We call the pair $(E,F)$ direction deficient if 
${\mathcal D}(E,F) \neq {\mathcal D}(\Bbb{F}_q^d)$. Thus there is a line through the origin $L$ 
such that $(E-F) \cap L \subseteq \{ \hat{0} \}$. 

To prove the first part of the theorem, we wish to show that any direction deficient pair 
$(E,F)$ has $|E| + |F| \leq q^d$ {\bf and} $|E \cap F| \leq q^{d-1}$. 
Once we show this it will follow that if 
$|E| + |F| > q^d$ {\bf or} $|E \cap F| > q^{d-1}$ 
then ${\mathcal D}(E,F)={\mathcal D}(\Bbb{F}_q^d)$ as desired.

For any linear change of variables given by a matrix $\mathbb{A} \in GL_d(\Bbb{F}_q)$, it is 
easy to see that $(E,F)$ is a direction deficient pair with ${\mathcal D}(E,F)$ missing line $L$ 
if and only if $(\mathbb{A}(E),\mathbb{A}(F))$ is a direction deficient pair with 
${\mathcal D}(\mathbb{A}(E), \mathbb{A}(F))$ missing line $\mathbb{A}(L)$.

Thus without loss of generality we can assume our direction deficient pair $(E,F)$ 
has $(E-F) \cap L \subseteq \{ \hat{0} \}$ where $L$ is the $x_d$-axis line i.e., the line through 
$(0,\dots,0,1)$.

Following the proof of Theorem \ref{graph}, we consider the projection 
$\pi: \Bbb{F}_q^d \to \Bbb{F}_q^{d-1}$ onto the first $(d-1)$ coordinates.
Let $\pi(E)=U$ and $\pi(F)=V$. Thus $U,V \subseteq \Bbb{F}_q^{d-1}$ however 
$\pi$ is not in general injective when restricted to $E$ or $F$ in this case.
Instead we will show that $\pi$ defines a bijection from $E \cap F$ to $U \cap V$.

Notice for $e \in E$ and $f \in F$ we have $\pi(e)=\pi(f)$ if and only if $e-f \in ker(\pi)=L$. 
However since $(E-F) \cap L \subseteq \{ \hat{0} \}$ we see this happens if and only if 
$e=f$ and so $e=f \in E \cap F$. From this it is easy to argue that 
$\pi(E \backslash F)=U \backslash V$, $\pi(F \backslash E)=V \backslash U$ and that $\pi |_{E \cap F} : E \cap F \to U \cap V$ is a bijection. 

Now since $|ker(\pi)|=|L|=q$ we have as $E \backslash F \subseteq \pi^{-1}(U \backslash V)$ 
that $|E \backslash F| \leq q |U \backslash V|$.
Similarly $|F \backslash E| \leq q |V \backslash U|$. 
Since $\pi: E \cap F \to U \cap V$ is a bijection, $|E \cap F| = |U \cap V|$.

Now we have
$$ |E| = |E \backslash F| + |E \cap F| \leq q|U \backslash V| + |U \cap V| $$
$$ |F| = |F \backslash E| + |E \cap F| \leq q|V \backslash U| + |U \cap V| $$
Thus as $q \geq 2$ we have 
$$ |E| + |F| \leq q( |U \backslash V| + |U \cap V| + |V \backslash U|)=q|U \cup V| \leq q(q^{d-1})=q^d $$
as $U \cup V$ is a subset of $\Bbb{F}_q^{d-1}$.
Furthermore $|E \cap F| = |U \cap V| \leq q^{d-1}$ as $U \cap V \subseteq \Bbb{F}_q^{d-1}$.

Thus we have shown that if $(E,F)$ is a direction deficient pair then 
$|E|+|F| \leq q^d$ {\bf and} $|E \cap F| \leq q^{d-1}$. Also note that in the proof, one can see that 
a direction deficient pair $(E,F)$, after a linear change of variables, has the intersection $E \cap F$ 
be a graph set of a function $f: U \cap V \to \Bbb{F}_q$ where $U,V \subseteq \Bbb{F}_q^{d-1}$ 
while $E \backslash F$ and 
$F \backslash E$ are contained in ruled solids (unions of lines parallel to $L$) 
lying above $U \backslash V$ and 
$V \backslash U$ respectively.

This completes the proof of the first part of the theorem as mentioned previously, it remains only to give an example that shows that the bound 
is sharp.

For this let $U, V$ be a partition of $\Bbb{F}_q^{d-1}$. Thus $U,V$ are disjoint and have union 
equal to $\Bbb{F}_q^{d-1}$. Then let $E=\pi^{-1}(U)$ and $F=\pi^{-1}(V)$ so 
$E, F$ is a partition of $\Bbb{F}_q^d$ and so $|E|+|F|=q^d$ and $|E \cap F|=0$. 
Note any element $e \in E$ is of the form $(u,x)$ for $u \in U, x \in \Bbb{F}_q$ and an element 
$f \in F$ is of the form $(v,y)$ for $v \in V, y \in \Bbb{F}_q$. Thus 
$e-f=(u-v,x-y)$ is never of the form $(0,\dots,0, \text{ nonzero })$ as $u \neq v$ since $U, V$ 
disjoint. Thus ${\mathcal D}(E,F)$ does not contain the $x_d$-axis line and 
so $(E,F)$ is a direction deficient pair with $|E|+|F|=q^d$.


\vskip.125in 

\begin{thebibliography}{10}

\vskip.125in

\bibitem{BKT04} J. Bourgain, N. Katz, and T. Tao, {\it A sum-product estimate in finite fields, and applications}, Geom. Funct. Anal. \textbf{14} (2004), 27-57.

\bibitem{D09} Z. Dvir, {\it On the size of Kakeya sets in finite fields}, J. Amer. Math. Soc. 22 (2009), no. 4, 1093-1097.

\bibitem{HLS10} D. Hart, L. Li, C. Y. Shen, {\it Fourier analysis and expanding phenomena in finite fields}, (http://arxiv.org/pdf/0909.5471). 

\bibitem{HIKR07} D. Hart, A. Iosevich, D. Koh and M. Rudnev, {\it Averages over hyperplanes, sum-product theory in vector spaces over finite fields and the Erdos-Falconer distance conjecture}, (accepted for publication by Transaction of the AMS), arXiv:0707.3473, (2007).

\bibitem{IR07} A. Iosevich and M. Rudnev, {\it Erd\H os distance problem in vector spaces over finite fields}, Trans. Amer. Math. Soc. \textbf{359} (2007), no. 12, 6127-6142.

\bibitem{LN97} R. Lidl and G. Niederreiter, {\it Finite fields}, Second edition, Cambridge University Press, (1997). 

\bibitem{PS04} J. Pach and M. Sharir {\it Geometric incidences}, Towards a theory of geometric graphs, 185-223, Contemp. Math., \textbf{342}, Amer. Math. Soc., Providence, RI, (2004). 

\bibitem{PPS07} J. Pach, R. Pinchasi, M. Sharir {\it Solution of Scott's problem on the number of directions determined by a point set in 3-space}, Discrete Comput. Geom. \textbf{38} (2007), no. 2, 399-441.

\bibitem{R99} I. Z. Ruzsa, {\it An analog of Freiman’s theorem in groups. Structure theory of set addition}, 
Aste ́risque \textbf{258} (1999), 323–326.

\bibitem{S08} T. Sanders, {\it A note on Freiman's theorem in vector spaces}, Combin. Probab. Comput. \textbf{17} (2)
(2008), 297–305.

\bibitem{TV06}  T. Tao and V. Vu, {Additive Combinatorics}, Cambridge University Press (2006).

\bibitem{V08} V. Vu, {\it Sum-product estimates via directed expanders}, Math. Res. Lett. \textbf{15} (2008), no. 2, 375-388.

\bibitem{W03} T. Wolff, {\it Lectures on harmonic analysis}. With a foreword by Charles Fefferman and preface by Izabella Laba. Edited by Laba and Carol Shubin. University Lecture Series, \textbf{29}, American Mathematical Society, Providence, RI, (2003).

\end{thebibliography}
\end{document}